\documentclass[11pt]{amsart}

\usepackage{amsmath}
\usepackage{amsxtra}
\usepackage{amscd}
\usepackage{amsthm}
\usepackage{amsfonts}
\usepackage{amssymb}
\usepackage{eucal}
\let\goth\mathfrak

\newcommand{\sll}{\widehat{\goth{sl}}_l}

\newcommand{\nn}{\nonumber}
\newcommand{\bea}{\begin{eqnarray}}
\newcommand{\ena}{\end{eqnarray}}
\newcommand{\be}{\begin{eqnarray*}}
\newcommand{\en}{\end{eqnarray*}}

\newcommand{\wt}{{\mathop{\rm wt}}}
\newcommand{\sgn}{{\mathop{\rm sgn}}}
\newcommand{\Ad}{{\mathop{\rm Ad}}}

\newcommand{\ch}{\mathop{{\rm ch}}}
\newcommand{\End}{\mathop{{\rm End}}}

\newcommand{\br}[1]{{\langle #1 \rangle}}  %bracket
\newcommand{\Wc}{\goth{S}_l}
\newcommand{\Pc}{\overset{\circ}{P}}
\newcommand{\Qc}{\overset{\circ}{Q}}

\newcommand{\W}{\mathcal{W}}
\newcommand{\Mc}{\mathcal{M}}

\newcommand{\bep}{\overline{\epsilon}}

\numberwithin{equation}{section}

\newtheorem{thm}{Theorem}[section]
\newtheorem{cor}[thm]{Corollary}
\newtheorem{prop}[thm]{Proposition}
\newtheorem{lem}[thm]{Lemma}
\theoremstyle{dfn}

\newtheorem{remark}[thm]{Remark}

\newcommand{\C}{{\mathbb C}}  
\newcommand{\Z}{{\mathbb Z}}

\begin{document}
\pagestyle{myheadings}
\markboth{Feigin, Jimbo, Loktev, Miwa and Mukhin}
{Monomial basis}

\title{Addendum to 
`Bosonic formulas for $(k,l)$-admissible partitions'}
\author{B.~Feigin, M.~Jimbo, S.~Loktev, T.~Miwa and E.~Mukhin}
\address{BF: Landau institute for Theoretical Physics, Chernogolovka,
142432, Russia}\email{feigin@feigin.mccme.ru}  
\address{MJ: Graduate School of Mathematical Sciences, University of
Tokyo, $\quad$ Tokyo 153-8914, Japan}\email{jimbomic@ms.u-tokyo.ac.jp}
\address{SL: Institute for Theoretical and Experemental Physics and Independent University of Moscow}\email{loktev@mccme.ru}
\address{TM: Division of Mathematics, Graduate School of Science, 
Kyoto University, Kyoto 606-8502
Japan}\email{tetsuji@kusm.kyoto-u.ac.jp}
\address{EM: Department of Mathematics, 
Indiana University-Purdue University-Indianapolis, 
402 N.Blackford St., LD 270, 
Indianapolis, IN 46202}
\email{mukhin@math.iupui.edu}

\date{\today}
\begin{abstract}
In our earlier paper we made a combinatorial study of 
$(k,l)$-admissible partitions. 
This object appeared already in the work of M. Primc 
as a label of a basis of level $k$-integrable modules over $\sll$. 
We clarify the relation between these two works. As a byproduct we
obtain an explicit parameterization of the affine Weyl group of $\sll$ by 
a simple combinatorial set.     
\end{abstract} 
\maketitle 

\medskip

\section{Introduction}\label{sec:1}
Let $\sll$ be the affine Lie algebra 
$\goth{sl}_l\otimes\C[t^{-1},t]\oplus\C c\oplus \C d$.  
In \cite{Pr}, M.~Primc 
constructed a basis of an integrable highest weight module $V(\Lambda)$
for any dominant integral weight $\Lambda$ of level $k$. 
With some change of notation, his basis is described as follows. 

We write $e_{ij}[n]=e_{ij}\otimes t^n$ where 
$e_{ij}=\left(\delta_{ia}\delta_{jb}\right)_{1\le a,b\le l}$. 
Consider the abelian Lie subalgebra $\goth{a}$ of $\sll$ spanned by 
$$
e_{21}[n],e_{31}[n],\cdots,e_{l1}[n], \qquad n\in\Z. 
$$
Let 
\be
W(\Lambda)=U(\goth{a})v_\Lambda 
\en
be the $\goth{a}$-submodule of $V(\Lambda)$ 
generated by the highest weight vector $v_\Lambda$, which satisfies
$e_{ij}[n]v_\Lambda=0$ ($n>0$).  
Denote by $\{\Lambda_i\}_{i=0}^{l-1}$ the fundamental weights and    
$\{\alpha_i\}_{i=0}^{l-1}$ the simple roots.
On $V(\Lambda)$ there is a projective representation of the lattice 
$\oplus_{i=1}^{l-1}\Z\alpha_i$, which we denote by 
$\gamma\mapsto T_\gamma$.  

Set 
\be
\beta=\alpha_{12}+\alpha_{13}+\cdots+\alpha_{1l},
%\label{beta} 
\en
where $\alpha_{ij}=\alpha_i+\alpha_{i+1}+\cdots+\alpha_{j-1}$ ($i<j$). 
Translating $W(\Lambda)$ by $T=T_\beta$,  
we have an increasing family of subspaces
$\cdots \subset T^mW(\Lambda)\subset T^{m+1}W(\Lambda)\subset\cdots$. 
It is not hard to see that the entire space is obtained as a limit:
\begin{prop}\label{prop:2.1}(\cite{Pr}, Theorem 8.2)
\be
V(\Lambda)=\br{T}W(\Lambda),
\en
where $\br{T}$ denotes the group generated by $T$. 
\end{prop}

A sequence ${\bf x}=(x_i)_{i=0}^\infty$ of 
integers with finitely many non-zero entries 
is called {\it $(k,l)$-admissible} if 
\be
0\le x_i\le k,
\quad 
x_i+x_{i+1}+\cdots+x_{i+l-1}\le k
\en
for all $i\ge 0$. 
Let $X(\Lambda)$ denote the set of all $(k,l)$-admissible sequences such that 
\be
&&
x_0+x_1+\cdots+x_i\le (\Lambda,\alpha_{1i+2}) 
\qquad (0\le i\le l-2).
%\label{adm2}
\en
For ${\bf x}\in X(\Lambda)$, introduce the vector 
\bea
M({\bf x})&=&
\prod_{n\ge 0}\prod_{i=2}^{l}e_{i1}[-n]^{x_{n(l-1)+i-2}}\,v_\Lambda
\label{mon}\\
&=&
\cdots 
e_{l1}[-1]^{x_{2l-3}}\cdots e_{21}[-1]^{x_{l-1}}
e_{l1}[0]^{x_{l-2}}\cdots e_{21}[0]^{x_0}v_\Lambda. 
\nn
\ena

\begin{thm}\label{thm:2.1}(\cite{Pr}, Theorem 9.1)
The set 
\bea
\Mc(\Lambda)=\{ M({\bf x}) \mid {\bf x}\in X(\Lambda)\}
\label{monbasis}
\ena
constitutes a basis of $W(\Lambda)$. 
\end{thm}
By Proposition \ref{prop:2.1},  
translating this basis one obtains a basis of $V(\Lambda)$. 
\begin{thm}\label{thm:2.2}(\cite{Pr}, Theorem 9.1) 
The space $V(\Lambda)$ has a basis
indexed by sequences $\tilde{{\bf x}}=(x_i)_{i\in\Z}$ 
of non-negative integers satisfying the following conditions:
\begin{enumerate}
\item $x_i=0$ for $i\gg 0$, 
\item $x_i=(\Lambda,\alpha_{i+1})$ for $i\ll 0$
(the index of $\alpha_{i+1}$ is to be read modulo $l$), 
\item $x_i+x_{i+1}+\cdots+x_{i+l-1}\le k$ for any $i\in\Z$.
\end{enumerate}
\end{thm}
Informally, one may think of $\tilde{{\bf x}}$ 
as representing a semi-infinite monomial
\be
&&\prod_{n\in \Z}\prod_{i=2}^{l}e_{i1}[-n]^{x_{n(l-1)+i-2}}\,v_{\infty}
\en
where $v_{\infty}$ is an ideal vector `at infinity'. 
Such a semi-infinite basis was also considered in \cite{FS}. 
Theorems \ref{thm:2.1} and \ref{thm:2.2} are the main results of \cite{Pr}.  

The basis \eqref{monbasis} is homogeneous with respect to 
the grading by the weight lattice 
$P=\oplus_{i=0}^{l-1}\Z\Lambda_i\oplus\Z\delta$
($\delta=\sum_{i=0}^{l-1}\alpha_i$).  
In general, for a $P$-graded vector space 
$U=\oplus_{\lambda\in P}U_\lambda$, we call the formal series 
\be
\chi (U)=\sum_{\lambda\in P}\dim (U_\lambda) e^{\lambda}
\en
the {\it character} of $U$. 
Theorem \ref{thm:2.1} implies that the character of $W(\Lambda)$ is given 
by the generating series of $(k,l)$-admissible partitions 
\bea
\chi\left(\Mc(\Lambda)\right)
=\sum_{{\bf x}\in X(\Lambda)}e^{\wt (M({\bf x}))},
\label{chi}
\ena
where the weight of the monomial \eqref{mon} is 
\be
\wt (M({\bf x}))=\Lambda-
\sum_{n\ge 0}\sum_{i=2}^{l}x_{n(l-1)+i-2}(\alpha_{1i}+n\delta).
\en
We studied the generating function \eqref{chi} 
and obtained a `bosonic' formula for it in \cite{FJLMM} 
(see also Remark \ref{rem:4.1} below). 
At the time of the writing, reference 
\cite{Pr} somehow escaped our notice. 
The purpose of this Addendum is to 
discuss the connection between the works \cite{Pr} and \cite{FJLMM}, 
and to present an alternative proof of Theorems 
\ref{thm:2.1} and \ref{thm:2.2}. 

The logic of our proof is as follows. 
Utilizing the relations arising from the integrability and 
highest weight condition,   
one shows that \eqref{monbasis} is a spanning set of
$W(\Lambda)$. This is an easy part of the proof 
which we do not repeat in this paper.

The spanning property together with Proposition \ref{prop:2.1} implies that 
\bea
\chi(V(\Lambda))
&=&\lim_{m\rightarrow\infty}\chi (T^mW(\Lambda))
\nn\\
&\ll&\lim_{m\rightarrow\infty}\chi\left(T^m \Mc(\Lambda)\right).
\label{cp2}
\ena
Here, 
for $f=\sum_{\lambda\in P}f_\lambda e^\lambda$ 
and $g=\sum_{\lambda\in P}g_\lambda e^\lambda$, 
we write $f\ll g$ if $f_\lambda\le g_\lambda$ holds for all $\lambda\in P$. 
Therefore the proof is completed if one shows that the 
last expression \eqref{cp2} coincides with $\chi( V(\Lambda))$. 
We do so by showing that in the limit \eqref{cp2}
the formula in \cite{FJLMM} gives rise to 
the Weyl-Kac character formula. This is the main result of this paper. 

It turns out that the limit can be taken
term-wise, no additional summation is needed.
In particular, each summand in the bosonic formula of \cite{FJLMM}
corresponds to a summand in the Weyl-Kac formula. In the former case,
the summands are parameterized by ``translated good monomials'' see
\eqref{good}, \eqref{translate good}, and in the latter case by the
elements of affine Weyl group of $\sll$. Therefore we obtain an explicit
bijection between these two sets, see Corollary \ref{bijection}.

\section{Bosonic formula}\label{sec:2}
\subsection{Notation}
Besides those introduced already, we use the following notation. 
Let $(~,~)$ be the invariant bilinear form 
on $\goth{h}^*=\oplus_{i=0}^{l-1}\C\Lambda_i\oplus\C\delta$
normalized as $(\alpha_i,\alpha_i)=2$.  
We use the orthonormal vectors $\epsilon_i$ ($i=1,\cdots,l$) 
in the Euclidean space 
to express $\overline{\Lambda}_i=\Lambda_i-\Lambda_0$ as 
$\overline{\Lambda}_i=\bep_1+\cdots+\bep_i$, 
where $\bep_i=\epsilon_i-\bep$, $\bep=\sum_{i=1}^l\epsilon_i/l$. 
We have $\alpha_{ij}=\bep_i-\bep_j$. 

The extended affine Weyl group $\widetilde{\W}\simeq \Pc\rtimes\Wc$    
is the subgroup of $GL(\goth{h}^*)$ 
generated by $\Pc=\oplus_{i=1}^{l-1}\Z\overline{\Lambda}_i$
and the symmetric group $\Wc$ on $l$ letters. 
The action of an element $w=\gamma\cdot \sigma\in\widetilde{\W}$
($\gamma\in \Pc$, $\sigma\in\Wc$) on $\goth{h}^*$ is given 
by the rule:
$w(\mu)=w\mu=t_\gamma(\sigma(\mu))$ ($\mu\in \goth{h}^*$),
where 
\bea
&&
\sigma(\Lambda_0)=\Lambda_0,~~ \sigma(\delta)=\delta, 
~~\sigma(\bep_i)=\bep_{\sigma(i)}\quad (1\le i\le l),\label{action1}
\\
&&
t_\gamma(\mu)=\mu+(\mu,\delta)\gamma -\left(\frac{1}{2}(\mu,\delta)
(\gamma,\gamma)+(\gamma,\mu)\right)\delta,  
\label{action2}
\ena
see formula (6.5.2) in \cite{K}.

The affine Weyl group $\W$ is the subgroup of the extended Weyl group 
$\W\simeq \Qc\rtimes\Wc$, $\Qc=\oplus_{i=1}^{l-1}\Z\alpha_i\subset
\Pc$. 

{}For a real root $\beta$, we denote by $r_\beta$ the corresponding reflection. Then we have
$r_\beta\in\W$, moreover, the affine Weyl group $\W$ 
is generated by reflections $r_{\alpha_i}$, $i=0,\dots,l-1$.
For $w_1,w_2\in\widetilde{\W}$ we set $\Ad_{w_1}(w_2)=w_1w_2w_1^{-1}$.
Then for $w\in\widetilde{\W}$ and a real root $\beta$ we have
\bea\label{ad on r}
\Ad_wr_\beta=r_{w(\beta)}.
\ena

\subsection{Previous results}\label{subsec:4.1}
Let us recall the main points of \cite{FJLMM}. 

Given a set of non-negative integers $\{b_i\}_{i=0}^{l-2}$ such that 
\be
0\le b_0\le b_1\le \cdots\le b_{l-2}\le k, 
\en
define a formal power series in the variables 
$q,z_1,\cdots,z_{l-1}$ 
\begin{eqnarray}
&&\chi_{b_0,\cdots,b_{l-2}}(q,z_1,\cdots,z_{l-1})
\nn\\
&&=\sum_{{\bf x}}z_1^{x_0}\cdots z_{l-1}^{x_{l-2}}
(qz_1)^{x_{l-1}}\cdots (qz_{l-1})^{x_{2l-3}}\cdots\,
\label{ch1}
\\
&&=\sum_{{\bf x}}\prod_{n\ge 0}\prod_{i=1}^{l-1}
(q^nz_i)^{x_{n(l-1)+i-1}},
\nn
\end{eqnarray}
where the sum runs over all 
$(k,l)$-admissible sequences ${\bf x}$
satisfying $x_0+x_1+\cdots+x_i\le b_i$
($0\le i\le l-2$). 
Eq. \eqref{ch1} is the unique power series solution of the difference equation 
\be
\chi_{b_0,\cdots,b_{l-2}}(q,z_1,\cdots,z_{l-1})
=
\sum_{i=0}^{b_0}z_1^{i} \chi_{b_1-i,\cdots,b_{l-2}-i,k-i}(q,z_2,\cdots,z_{l-1},qz_1)
\en
with the initial condition 
$\chi_{b_0,\cdots,b_{l-2}}(q,0,\cdots,0)=1$. 

Let $P_1,\cdots,P_l$ be monomials in 
$q^{\pm 1},z_1^{\pm 1},\cdots,z_{l-1}^{\pm 1}$, 
and write 
\be
[P_1,\cdots,P_l]=P_1^{b_0}P_2^{b_1-b_0}\cdots P_l^{k-b_{l-2}}. 
\en

We define linear operators $A,B$ which act
on expressions of the form $f(q,z_1,\cdots,z_{l-1})[P_1,\cdots,P_l]$
by the formulas
\begin{align*}
A\left(f[P_1,\cdots,P_l]\right)
&=\frac{S(f)}{1-z_1S(P_l/P_1)}[S(P_1),S(P_1),S(P_2),\cdots,S(P_{l-1})],
\\
B\left(f[P_1,\cdots,P_l]\right)
&=\frac{S(f)}{1-z_1^{-1}S(P_1/P_l)}[z_1S(P_l),S(P_1),S(P_2),\cdots,S(P_{l-1})],
\end{align*}
where $S$ stands for the shift operator
\bea
(Sf)(q,z_1,\cdots,z_{l-1})=f(q,z_2,\cdots,z_{l-1},qz_1).
\label{shift}
\ena
In terms of $A,B$, 
the difference equation and the initial condition 
can be recast into the form 
\be
\chi_{b_0,\cdots,b_{l-2}}
=\lim_{m\rightarrow\infty}(A+B)^m [1,\cdots,1].
\en

An ordered product of operators $A,B$,  
$M=C_1C_2\cdots C_n$ ($C_i\in\{A,B\}$),  
is called a monomial, and $n$ is called the degree of $M$.  
When $n=0$, $M=1$ means the identity.  
$M$ is called {\em good} if 
\bea\label{good}
\mbox{$C_i=A$ implies $C_{i+l-1}=A$ for $i=1,\cdots,n-l+1$.}
\ena 
Let $\mathcal{G}$ denote the set of good monomials such that 
either $n=0$, or else $C_n=B$.  
The main result of  \cite{FJLMM} is the formula  
\bea
\chi_{b_0,\cdots,b_{l-2}}=
\sum_{M\in\mathcal{G}}M{\bf v}_\infty,
\label{bos1}
\ena
where
\be
{\bf v}_{\infty}=\frac{1}{(z_1)_\infty\cdots(z_{l-1})_\infty}[1,\cdots,1],
\en
and $(z)_\infty=\prod_{n=0}^\infty(1-q^nz)$. 

\begin{remark}\label{rem:4.1}
Actually the character treated in \cite{FJLMM} is 
a specialization of \eqref{ch1} given by 
$$
q=\overline{q}^{l-1},~
z_1=z,~
z_2=\overline{q}z,\cdots,
z_{l-1}=\overline{q}^{l-2}z,
$$
where $\overline{q},z$ stand for the variables $q,z$ used in \cite{FJLMM}. 
The working of \cite{FJLMM} carries over straightforwardly 
to the present setting as well. 
\end{remark}

\subsection{Reformulation by root systems}\label{subsec:4.2}
Let us reformulate the above results in terms of the root system of $\sll$. 

In the sequel we make the identification of variables
\be
q=e^{-\delta},z_1=e^{-\alpha_{12}},\cdots,z_{l-1}=e^{-\alpha_{1l}}.
%\label{qz}
\en
Then a monomial in 
$q^{\pm 1},z_1^{\pm 1},\cdots,z_{l-1}^{\pm 1}$ 
is an element of the group ring $\C[Q]$. 
The shift operator \eqref{shift} acting on $\C[Q]$ is 
implemented as $S(e^\beta)=e^{s(\beta)}$,  
where $s$ is the element of the extended Weyl group given by
\be
s=\bep_2\cdot \sigma,
\quad \sigma=(2\,3\,\cdots\,l)\in\goth{S}_l.
\en

{}To each good monomial $M\in\mathcal{G}$, we associate 
an element $w_M$ of the affine Weyl group by the following rule.
Write $M$ as a product $C_1C_2\cdots C_n$ and read it from the left. 
Associate a sequence of roots 
$\beta_1,\beta_2\cdots=\alpha_{12},\alpha_{13},\cdots$
with 
\be
\beta_{m(l-1)+r}=\alpha_{1\,r+1}+m\delta
\qquad (m\in\Z, 1\le r\le l-1).
%\label{beta}
\en
{}From the above sequence, drop $\beta_i$ such that $C_i=A$, and 
denote the resulting sequence by $\gamma_1,\cdots, \gamma_p$. 
We define 
\be
w_M=r_{\gamma_p}\cdots r_{\gamma_1}.
%\label{gammap} 
\en
{}For example, if $l=3$ and $M=BBBABAB$, then 
\be
w_M=
(r_{\alpha_{12}+3\delta})
(r_{\alpha_{12}+2\delta})
(r_{\alpha_{12}+\delta})
r_{\alpha_{13}} 
r_{\alpha_{12}}.
\en
Alternatively, $w_M$ is determined recursively as follows. 
\begin{lem}\label{lem:2.5}
\be
&&w_{1}=1,
\\
&&
w_{AM}=\Ad_s(w_M),
\\
&&
w_{BM}=\Ad_s(w_M)r_{\alpha_{12}}.
\en
\end{lem}
\begin{proof}
This follows from the definition, formula \eqref{ad on r}
and 
\be
&&s(\alpha_{12})=\alpha_{13},
~~s(\alpha_{13})=\alpha_{14},\cdots,
s(\alpha_{1l-1})=\alpha_{1l},\quad 
s(\alpha_{1l})=\alpha_{12}+\delta,
\\
&& s(\bep_1)=\bep_1+\frac{1}{l}\delta,
\qquad s(\delta)=\delta,
\en
obtained from formulas \eqref{action1}, \eqref{action2}.
\end{proof}

Next we define for each $M\in\mathcal{G}$ 
a formal Laurent series $f_M$ in $q,z_1,\cdots,z_{l-1}$ 
inductively as follows.         
\bea
&&
f_1=\frac{1}{(e^{-\alpha_{12}})_\infty\cdots (e^{-\alpha_{1l}})_\infty},
\label{scalar0}
\\
&&f_{BM}=
S(f_{M})\frac{1}{1-e^{\Ad_s(w_M)\bep_1-\bep_2}},
\label{scalar1}
\\
&&
f_{AM}=S(f_{M})\frac{1}{1-e^{-\Ad_s(w_M)\alpha_{12}}}
=-e^{\Ad_s(w_M)\alpha_{12}}f_{BM}.
\label{scalar2}
\ena
We apply the last formula when $AM$ is also a good monomial. 
In \eqref{scalar2} we have used 
\bea
\Ad_s(w_M)\bep_2=\bep_2~~~ \mbox{ if }~~~ AM\in\mathcal{G}. 
\label{ads} 
\ena
Formula \eqref{ads} follows from the fact that if $AM\in\mathcal{G}$ 
a root of the form $\alpha_{12}+j\delta$ does not appear 
in the sequence used to calculate $\Ad_s(w_M)$.  

The result \eqref{bos1} is then rephrased as follows. 
\begin{prop}\label{prop:4.1}
{}For a dominant integral weight $\Lambda$ of level $k=(\Lambda,\delta)$, 
we have 
\be
\chi\left(\Mc(\Lambda)\right)=\sum_{M\in\mathcal{G}}f_M\,e^{w_M(\Lambda)}.
\en
\end{prop}

\begin{proof}
We prove 
\bea
M{\bf v}_\infty=f_{M}\,e^{w_M(\Lambda)-\Lambda} 
\label{M}
\ena
by induction on the degree of $M$. 
The Proposition is an immediate consequence of \eqref{bos1} and 
\eqref{M}. 

In general, for $\lambda\in P$ and $w\in \W$ we have
\be
w(\lambda)-\lambda=\sum_{i=0}^{l-1}(\lambda,\alpha_i)(w(\Lambda_i)-\Lambda_i).
\en
Upon writing $\gamma_i=w(\Lambda_i)-\Lambda_i$ and noting 
$\Ad_s(w)\lambda-\lambda=s(ws^{-1}\lambda-s^{-1}\lambda)$, 
it follows that 
\bea
\Ad_s(w)\lambda-\lambda&=&
(\lambda,\alpha_1+\alpha_2)s\gamma_1
+\sum_{i=2}^{l-2}(\lambda,\alpha_{i+1})s\gamma_i
\nn\\
&&+(\lambda,\alpha_0+\alpha_1)s\gamma_{l-1}-(\lambda,\alpha_1)s\gamma_0,
\label{adsw1}\\
\Ad_s(w)r_{\alpha_1}\lambda-\lambda&=&
(\lambda,\alpha_1)s(\alpha_0+\gamma_0)+
\sum_{i=1}^{l-2}(\lambda,\alpha_{i+1})s\gamma_i
\nn\\
&&+(\lambda,\alpha_0)s\gamma_{l-1}.
\label{adsw2}
\ena

Take $b_i=(\Lambda,\alpha_{1i+2})$ ($0\le i\le l-2$), so that 
\be
\Lambda=(k-b_{l-2})\Lambda_0+b_0\Lambda_1+(b_1-b_0)\Lambda_2+
\cdots+(b_{l-2}-b_{l-3})\Lambda_{l-1}.
\en 
{}For a good monomial $M\in\mathcal{G}$ we set 
$P_i=e^{\gamma_i}$, where 
$\gamma_i=w_M(\Lambda_i)-\Lambda_i$ and $P_l=P_0$. 
In the notation of the previous subsection and using 
\eqref{adsw2},   
we have then 
\be
&&
[P_1,\cdots,P_l]=e^{w_M(\Lambda)-\Lambda},
\\
&&
[z_{1}S(P_l),S(P_1),S(P_2),\cdots,S(P_{l-1})]=
e^{w_{BM}(\Lambda)-\Lambda},
\\
&&
f_{BM}
=\frac{S(f_M)}{1-z_{1}^{-1}S(P_1/P_l)}.
\en
In the case $AM$ is also a good monomial, we find from \eqref{ads} that 
$w_M(\Lambda_{l-1})-\Lambda_{l-1}=w_M(\Lambda_0)-\Lambda_0$. 
Eq. \eqref{adsw1} then simplifies to yield 
\be
&&
[S(P_1),S(P_1),S(P_2),\cdots,S(P_{l-1})]=e^{w_{AM}(\Lambda)-\Lambda},
\\
&&
f_{AM}=\frac{S(f_M)}{1-z_{1}S(P_l/P_1)}. 
\en

Therefore
\be
f_{BM}e^{w_{BM}(\Lambda)-\Lambda}&=&
B\left(f_Me^{w_M(\Lambda)-\Lambda}\right),
\\
f_{AM}e^{w_{AM}(\Lambda)-\Lambda}&=&
A\left(f_Me^{w_M(\Lambda)-\Lambda}\right).
\en
\end{proof}

Introduce the formal inverse ofi letter $B$, $B^{-1}$, 
with the property $B^{-1}B=BB^{-1}=1$ 
and consider the set of {\em translated good monomials}
\bea\label{translate good}
\widetilde{\mathcal{G}}=\{B^{-n}M\mid n\ge 0,~M\in\mathcal{G}\}.
\ena

We extend the definition of $w_M$, $f_M$ for 
$M\in\widetilde{\mathcal{G}}$ by the rule
\be
&&w_{B^{-1}M}=(\Ad_s)^{-1}(w_Mr_{\alpha_{12}}),
\\
&&f_{B^{-1}M}=S^{-1}\left((1-e^{w_M(\bep_2)-\bep_2})f_M\right).
\en
Note that $w_{B^{-1}BM}=w_{BB^{-1}M}=w_M$ and 
$f_{B^{-1}BM}=f_{BB^{-1}M}=f_M$, so $w_M$ and
$f_M$, $M\in\widetilde{\mathcal{G}}$, are well defined.

\begin{lem}\label{lem:3.2}
The map $w:\;\widetilde{\mathcal{G}}\rightarrow \W$ sending 
$M$ to $w_M$ is injective. 
\end{lem}
\begin{proof}
The restriction of the map $w$ to $\mathcal{G}$ is injective.
Indeed, let the sequence 
${\bf x}=(x_i)_{i=0}^\infty$ 
be the extremal configuration in the sense of 
\cite{FJLMM}, Subsection 2.5. Then we have 
\be
e^{w_M(\Lambda)-\Lambda}=\prod_{n=0}^\infty \prod_{i=1}^{l-1}(q^nz_i)^{x_{n(l-1)+i-1}}.
\en
It is proved in Proposition 2.11 there that 
if we regard $b_i=(\Lambda,\alpha_{1i+2})$ as variables, then 
a good monomial $M\in \mathcal{G}$ is uniquely determined from ${\bf x}$. 
Since $\Lambda$ is arbitrary, our claim follows.

Notice that if $w_{M_1}\neq w_{M_2}$ then $w_{B^{-1}M_1}\neq
w_{B^{-1}M_2}$. Now the lemma follows from the fact that if
$M\in\mathcal{G}$ then $B^nM\in\mathcal{G}$ for $n>0$.
\end{proof}

We set 
\be
r^{(j)}=
\begin{cases} 
r_{s^{j-1}(\alpha_{12})}\cdots r_{s(\alpha_{12})}r_{\alpha_{12}} & (j\ge 0), \\
r_{s^j(\alpha_{12})}\cdots r_{s^{-2}(\alpha_{12})}r_{s^{-1}(\alpha_{12})} & (j< 0), \\
\end{cases}
%\label{betaseq}
\en
so that $r^{(i+j)}=(\Ad_s)^i(r^{(j)})r^{(i)}$ and 
$w_{B^jM}=(\Ad_s)^j(w_M)r^{(j)}$ for all $i,j\in \Z$.  

\begin{lem}\label{lem:ext}
\bea
s^{l(l-1)}=r^{(l(l-1))}=t_{-\beta}, 
\label{srt}
\ena
where $\beta=\alpha_{12}+\cdots+\alpha_{1l}=l\bep_1$. 
\end{lem}
\begin{proof}
This follows from the relations
\be
&&s^j=t_{\bep_2+\cdots+\bep_{j+1}}\cdot\sigma^j\qquad (1\le j\le l-1),\\
&&r_{\alpha_{1j}-n\delta}=t_{n\alpha_{1j}}r_{\alpha_{1j}},\\
&&r^{(j(l-1))}=
t_{-\alpha_{12}-\cdots-\alpha_{1\,j}}\cdot(12\cdots l)^{j}
\qquad (1\le j\le l).
\en
\end{proof}

\subsection{Comparison with the Weyl-Kac formula}\label{subsec:4.3}
Recall the translation operator $T\in\End\left(V(\Lambda)\right)$  
which gives a linear isomorphism between weight spaces 
$T:V(\Lambda)_\lambda\rightarrow V(\Lambda)_{t_\beta\lambda}$.  
We note that acting with the operator $B^l$ corresponds to $T$. 
If we write     
\be
\chi\left(\Mc(\Lambda)\right)=\sum_{{\bf x}\in X(\Lambda)}e^{\wt (M({\bf x}))}
=e^\Lambda\chi(q,z_1,\cdots,z_{l-1}),
\en
then the character of the translated set $T^m\Mc(\Lambda)$ becomes 
\be
\chi\left(T^m\Mc(\Lambda)\right)
&=&\sum_{{\bf x}\in X(\Lambda)}e^{t_{m\beta}(\wt (M({\bf x})))}
\\
&=&
e^{t_{m\beta}(\Lambda)}\chi(q,q^{-ml}z_1,\cdots,q^{-ml}z_{l-1})
\\
&=&
e^{t_{m\beta}(\Lambda)}(S^{-ml(l-1)}\chi)(q,z_1,\cdots,z_{l-1}).
\en
{}From Proposition \ref{prop:4.1}, we obtain 
\be
\chi\left(T^m \Mc(\Lambda)\right)
&=&
e^{t_{m\beta}(\Lambda)}
\sum_{M\in\mathcal{G}}S^{-ml(l-1)}\left(f_Me^{w_M(\Lambda)-\Lambda}\right)
\\
&=&
\sum_{M\in\mathcal{G}}S^{-ml(l-1)}(f_M)
e^{(t_{m\beta} w_M)(\Lambda)}
\\
&=&
\sum_{M\in B^{-ml(l-1)}\mathcal{G}}
S^{-ml(l-1)}(f_{B^{ml(l-1)}M})e^{w_M(\Lambda)},
\en
where we have used Lemma \ref{lem:ext}.
If we set 
\be
h_M=\lim_{n\rightarrow\infty}S^{-n}(f_{B^nM}),
\en
then
\bea
\lim_{m\rightarrow\infty}\chi\left(T^m\Mc(\Lambda)\right)
=\sum_{M\in\widetilde{\mathcal{G}}}h_M\cdot e^{w_M(\Lambda)}.
\label{lim1}
\ena

\begin{lem}\label{lem:4.1}
\bea
\frac{f_M}{h_M}&=&
\prod_{2\le i\le l}(qe^{w_M(\bep_1)-\bep_i})_\infty
\label{product}
\\
&\times&
\prod_{2\le j<i\le l}(e^{w_M(\bep_i)-\bep_j})_\infty
\prod_{2\le i\le j\le l}(q e^{w_M(\bep_i)-\bep_j})_\infty.
\nn
\ena
\end{lem}

\begin{proof}
Applying \eqref{scalar1} repeatedly we obtain 
\be
S^{-n}(f_{B^nM})
=f_M
\prod_{j=1}^n
(1-e^{(w_M s^{-j}r^{(j)})(\bep_2)-s^{-j}(\bep_2)})^{-1}.
\en
To calculate $s^{-j}(\bep_2)$ 
note relations \eqref{srt} and 
\be
&&
s^{-(l-1)j-p}(\bep_2)=\bep_{l+1-p}+\frac{l+j-p}{l}\delta,
\quad (j\in\Z, 0\le p\le l-2).  
\en
{}For the calculation of $(s^{-j}r^{(j)})(\bep_2)$ 
we use 
\be
r^{((l-1)j+p)}(\bep_2)
=\begin{cases}
\bep_{j+3}& (0\le j\le p-2),\\
\bep_{1}+j\delta& (j=p-1),\\
\bep_{j+2}& (p\le j\le l-2),\\
\bep_{2}-\delta& (j=l-1),\\
\end{cases}
\en
where $1\le p\le l-1$. 
After simplification 
we obtain formula \eqref{product}.
\end{proof}

Let $D$ denote the Weyl denominator
\be
&&D=(q)_\infty^{l-1}
\prod_{i=1}^{l-1}(z_i)_\infty(qz_i^{-1})_\infty
\prod_{1\le i<j\le l-1}
(z_i^{-1}z_j)_\infty(qz_iz_j^{-1})_\infty. 
\en

Let $\rho=\Lambda_0+\cdots+\Lambda_{l-1}$.
For all $w\in \W$, we have  
\bea
w(e^\rho D)=\sgn (w) \,e^\rho D.
\label{D}
\ena

\begin{lem}\label{lem:4.2}
{}For $M\in\widetilde{\mathcal{G}}$ we have 
\be
h_M=w_M(D^{-1})=D^{-1}{\rm sgn}\, (w_M) \,e^{w_M(\rho)-\rho}.
\en
\end{lem}
\begin{proof}
The second equality is just a version of \eqref{D}.
{}From the definition of $h_M$ it is clear that 
\be
&&h_{BM}=S(h_M).
\en
Using Lemma \ref{lem:4.1}, we obtain 
\be
&&h_1=D^{-1},
\\
&&h_{AM}=-h_{BM}e^{\Ad_s(w_M)\alpha_{12}}.
\en
These properties uniquely characterize $h_M$. 
On the other hand, a direct verification shows 
that $S(D)=-D e^{\alpha_{12}}$. 
Noting $s(\rho)=\rho-\alpha_{12}$, we have for all $w\in\W$
\be
&&\Ad_s(w)r_{\alpha_{12}}(D^{-1})=S(w(D^{-1})),
\\
&&\Ad_s(w)(D^{-1})=-S(w(D^{-1}))e^{\Ad_s(w)\alpha_{12}}.
\en
The lemma follows from these equalities. 
\end{proof}

\begin{prop}\label{prop:4.2}
\be
\lim_{m\rightarrow\infty}\chi\left(T^m\Mc(\Lambda)\right)
=\ch (V(\Lambda)).
\en
\end{prop}

\begin{proof}
{}From \eqref{lim1} and Lemma \ref{lem:4.2}, we find 
\be
\lim_{m\rightarrow\infty}\chi\left(T^m\Mc(\Lambda)\right)
=D^{-1}
\sum_{M\in\widetilde{\mathcal{G}}}\sgn (w_M)
e^{w_M(\Lambda+\rho)-\rho}.
\en
By Lemma \ref{lem:3.2}, 
the right hand side can be written as a sum 
$\sum_{w\in\mathcal{S}}$ over some subset 
$\mathcal{S}$ of $\W$. 
In the special case $k=0$, the left hand side is $1$, 
so the Weyl denominator formula implies
\be
\sum_{w\in\mathcal{S}}\mathop{\rm sgn}(w)\, e^{w(\rho)-\rho}
=\sum_{w\in\W}\mathop{\rm sgn}(w)\, e^{w(\rho)-\rho}.
\en
Since $w(\rho)=\rho$ holds if and only if $w=1$, 
we must have that $\mathcal{S}=\W$. 
The proposition is then a consequence
of the Weyl-Kac character formula (\cite{K}, formula 10.4.1).
\end{proof}

\begin{cor}\label{bijection}
The map $w:\;\widetilde{\mathcal{G}}\rightarrow \W$ sending 
$M$ to $w_M$ is bijective.
\end{cor}

\bigskip

%%%%%%%%%%%%%%%%%%%%%%%%%%%%%%%%%%%%%%%%%%%%%%%%%%%%%%%%%%%%%%%
%                                                             %
%  Acknowledgments                                            %
%                                                             %
%%%%%%%%%%%%%%%%%%%%%%%%%%%%%%%%%%%%%%%%%%%%%%%%%%%%%%%%%%%%%%%
\noindent
{\it Acknowledgments.}\quad 
This work is partially supported by
the Grant-in-Aid for Scientific Research (B)
no.12440039 and (A1) no.13304010, Japan Society for the Promotion of Science.
The work of SL is partially supported by the grant RFBR-01-01-00546.

%%%%%%%%%%%%%%%%%%%%%%%%%%%%%%%%%%%%%%%%%%%%%%%%%%%%%%%%%%%%%%%
%                                                             %
%  References                                                 %
%                                                             %
%%%%%%%%%%%%%%%%%%%%%%%%%%%%%%%%%%%%%%%%%%%%%%%%%%%%%%%%%%%%%%%
%\bibliographystyle{unsrt}
%\bibliography{q}

\end{document}